\documentclass[10pt,a4paper,twoside]{article}

\usepackage{amssymb,amsmath,multirow}

\newtheorem{theorem}{Theorem}[section]

\newtheorem{lemma}[theorem]{Lemma}

\newtheorem{remark}[theorem]{Remark}

    \makeatletter
    \let\@fnsymbol\@arabic
    \makeatother

\title{Three Characterizations of Exponential Distribution Involving Median of Sample of Size Three}
\author{Marko Obradovi\' c}
\date{}

\begin{document}
\maketitle
\begin{abstract}
In this paper three new characterizing theorems of exponential
distribution are presented. They are based on equidistribution of some functions of order statistics. All of them include the median of sample of size three.
\end{abstract}

{\small \textbf{ keywords:} characterization, exponential distribution, order statistics

}

\section{Introduction}

The characterization theorems have been very popular lately. Since exponential distribution has  very wide applications, a major part of work on characterizations has been dedicated to this distribution.
Many examples could be found in \cite{ahsanullah} and \cite{galambos}.

Recently Arnold and Villasenor \cite{arnold} proposed a series of characterizations based on sample of size  two and made some conjectures on their generalization. 
They also proposed a new method of proof which can be used when the density in question is analytic. Later Yanev and Chakraborty \cite{yanev}
proved by this method two characterization theorems concerning maximum of sample of size three.

One of the benefits of these characterizations for small sample sizes is that they are suitable for creation of goodness of fit tests. 

We propose three new characterization theorems whose common point is that in all of them appears the median of sample of size three. In some sense these characterizations could be considered as a generalization
of the characterizations from \cite{arnold} and \cite{yanev}.

\section{Main results}

 Let $\mathcal{F}$ be a class of absolutely continuous distribution functions $F$ such that $F(0)=0$ and whose density function $f$ allows expansion in Maclaurin's series for all $x>0$.

 Firstly we  present two lemmas which are important for the proofs of the theorems.

\begin{lemma}\label{exp}
 Let $F$ be a distribution function from $\mathcal{F}$. If for all natural $k$ holds:

 \begin{equation}\label{fizvod}
  f^{(k)}(0)=(-1)^kf^{k+1}(0),
 \end{equation}
then $f(x)=\lambda e^{-\lambda x}$ for some $\lambda>0$.
\end{lemma}
\textbf{Proof.} Expanding the function $f$ in Maclaurin series for positive values of $x$ we get:
\begin{equation}
 f(x)=\sum\limits_{k=0}^{\infty}f^{(k)}(0)\frac{x^{k}}{k!}=\sum\limits_{k=0}^{\infty}(-1)^kf^{k+1}(0)\frac{x^{k}}{k!}=f(0)e^{-f(0)x}.
\end{equation}
For $f(0)>0$ this is the density of exponential distribution with $\lambda=f(0)$. $\hfill \Box$

\begin{lemma}\label{GiH}
 Let $F$ be a distribution from $\mathcal{F}$. Denote $G(x)=F(x)f(x)$ and $H(x)=F^2(x)f(x)$.
 Let the condition \eqref{fizvod}
 be satisfied for $k\leq r-2$. Then the following equalities hold:
 \begin{eqnarray}
  G^{(k)}(0)&=&(-1)^{k-1}f^{k+1}(0)(2^k-1), \;\;1\leq k\leq r-1.\\
  H^{(k)}(0)&=&(-1)^{k-2}f^{k+1}(0)(3^k-2^{k+1}+1), \;\;1 \leq k\leq r.\label{Hx}
  \end{eqnarray}

\end{lemma}
\textbf{Proof.} Applying Leibniz formula for the derivative of a product to $G(x)$ we get:

\begin{equation}
 G^{(k)}(x)=\sum\limits_{j=0}^{k}\binom{k}{j}F^{(j)}(x)f^{(k-j)}(x).
\end{equation}

When we put $x=0$, all summands that have $F^{(0)}(x)$ will be equal to zero, therefore we get

\begin{eqnarray*}
 G^{(k)}(0)&=&\sum\limits_{j=1}^{k}\binom{k}{j}f^{(j-1)}(0)f^{(k-j)}(0)\\
            &=&\sum\limits_{j=1}^{k}\binom{k}{j}(-1)^{j-1}f^{j}(0)(-1)^{k-j}f^{k-j+1}(0)\\
            &=&(-1)^{k-1}f^{k+1}(0)\sum\limits_{j=1}^{k}\binom{k}{j}\\
            &=&(-1)^{k-1}f^{k+1}(0)(2^k-1).
\end{eqnarray*}

Note that the expression on the right hand side of \eqref{Hx} is equal to zero for $k=1$.
Since $H(x)=F^2(x)f(x)=F(x)F(x)f(x)$, applying Leibniz formula for the derivative of product of three functions we obtain:
\begin{equation}
 H^{(k)}(x)=\sum\limits_{s=0}^{k}\sum\limits_{j=0}^{s}\binom{k}{j,s-j,k-s}F^{(j)}(x)F^{(s-j)}(x)f^{(k-s)}(x).
\end{equation}

Putting $x=0$ and using the same argument as in the previous case we get that $H'(0)=0$ and for $k\geq 2$

\begin{eqnarray*}
H^{(k)}(0)\!\!\!\!\!&=&\!\!\!\!\!\sum\limits_{s=2}^{k}\sum\limits_{j=1}^{s}\binom{k}{j,s-j,k-s}f^{(j-1)}(0)f^{(s-j-1)}(0)f^{(k-s)}(0)\\
           &=&\!\!\!\!\!\sum\limits_{s=2}^{k}\sum\limits_{j=1}^{s-1}\binom{k}{j,s-j,k-s}(-1)^{j-1}f^{j}(0)(-1)^{s-j-1}f^{s-j}(0)(-1)^{k-s}f^{k-s+1}(0)\\
           &=&\!\!\!\!\!(-1)^{k-2}f^{k+1}(0)\sum\limits_{s=2}^{k}\sum\limits_{j=1}^{s-1}\binom{k}{j,s-j,k-s}\\
           &=&\!\!\!\!\!(-1)^{k-2}f^{k+1}(0)\left(3^k-\sum\limits_{s=0}^{1}\sum\limits_{j=0}^{s}\binom{k}{j,s-j,k-s}-2\sum\limits_{s=2}^{k}\binom{k}{s}\right)\\
           &=&\!\!\!\!\!(-1)^{k-2}f^{k+1}(0)(3^k-2^{k+1}+1).
\end{eqnarray*}

$\hfill \Box$

\begin{remark}
 Note that the expressions for $G'(0)=f^{2}(0)$ and $H''(0)=2f^{3}(0)$ are always true regardless of the condition of the theorem. Besides, $G(0)=H(0)=H'(0)=0$.
\end{remark}

Let $X_{(k;n)}$ be the $k$ th order statistics from the sample of size $n$.
We now state and prove the characterization theorems.

\begin{theorem}\label{prva}
Let $X_1,X_2,X_3$ be independent random variables from the distribution $F\in\mathcal{F}$. If

\begin{equation}
 \frac{1}{3}X_1+\frac{1}{2}X_2\overset{d}{=} X_{(2;3)}
\end{equation}

then $X\sim \mathcal{E}(\lambda), \lambda>0$.

\end{theorem}
\textbf{Proof.}
Equalizing the densities on the left and the right hand side we get
\begin{eqnarray}
\nonumber&&\int\limits_{0}^x3f(3y)2f(2(x-y))dy=6F(x)(1-F(x))f(x)\\
\nonumber&&\int\limits_{0}^xf(3y)f(2(x-y))dy=f(x)\int\limits_{0}^xf(y)(1-2F(y))dy\\
&&\int\limits_{0}^xf(3y)f(2(x-y))dy=f(x)\int\limits_{0}^x(f(y)-2G(y))dy.\label{t1indegralna}
\end{eqnarray}

We shall prove by induction that from \eqref{t1indegralna} follows \eqref{fizvod}.
Differentiating both sides twice with respect to $x$ we get
\begin{eqnarray*}
&&f'(3x)3f(0)+2f(3x)f'(0)+4\int\limits_{0}^xf(3y)f''(2(x-y))dy\\&=&2f(x)f'(x)-2G'(x)f(x)-2G(x)f'(x)+(f(x)-2G(x))f'(x)\\&+&f''(x)\int\limits_{0}^{x}(f(y)-2G(y))dy.
\end{eqnarray*} 
Letting $x=0$ we get 
\begin{equation*}
5f'(0)f(0)=3f(0)f'(0)-2f(0)G'(0).
\end{equation*}
Using lemma \ref{GiH} we get 
\begin{equation*}
f'(0)=-f^2(0).
\end{equation*}
Hence \eqref{fizvod} is satisfied for $k=1$.
Let us now suppose that \eqref{fizvod} is satisfied for $1\leq k\leq r-2$. Differentiating both sides of \eqref{t1indegralna} $r$ times we get
\begin{eqnarray*}
&&\sum\limits_{j=0}^{r-1}3^{r-1-j}f^{(r-1-j)}(3x)2^jf^{(j)}(0)+\int\limits_{0}^xf(3y)2^rf^{(r)}(2(x-y))dy\\
&=&\sum\limits_{j=1}^r\binom{r}{j}f^{(r-j)}(x)(f^{(j-1)}(x)-2G^{(j-1)}(x))+f^{(r)}(x)\int\limits_{0}^x(f(y)-2G(y))dy.
\end{eqnarray*}
Letting $x=0$ and eliminating zero terms we get 
\begin{eqnarray}
\sum\limits_{j=1}^{r}3^{r-j}f^{(r-j)}(0)2^{j-1}f^{(j-1)}(0)=\sum\limits_{j=1}^r\binom{r}{j}f^{(r-j)}(0)(f^{(j-1)}(0)-2G^{(j-1)}(0)).
\end{eqnarray}
Extracting the summands which contain $f^{(r-1)}(0)$ and grouping them to one side we get
\begin{eqnarray*}
&&f^{(r-1)}(0)f(0)(r+1-3^{r-1}-2^{r-1})\\&=&\sum\limits_{j=2}^{r-1}f^{(r-j)}(0)3^{r-j}2^{j-1}f^{(j-1)}(0)-\sum\limits_{j=2}^{r-1}\binom{r}{j}
f^{(r-j)}(0)(f^{(j-1)}(0)-2G^{(j-1)}(0)).
\end{eqnarray*}
Using the induction hypothesis and lemmas \ref{exp} and \ref{GiH} we obtain
\begin{equation*}
f^{(r-1)}(0)(r+1-3^{r-1}-2^{r-1})=(-1)^{r-1}f^r(0)\bigg(\sum\limits_{j=2}^{r-1}\Big(3^{r-j}2^{j-1}-\binom{r}{j}2^{j-1}\Big)-2(2^{r-1}-1)\bigg).
\end{equation*}
To prove that \eqref{fizvod} is satisfied for $k=r-1$ it suffices to show that 
 \begin{equation*}
r+1-3^{r-1}-2^{r-1}=\sum\limits_{j=2}^{r-1}\Big(3^{r-j}2^{j-1}-\binom{r}{j}2^{j-1}\Big)-2(2^{r-1}-1),
\end{equation*}
or equivalently 
\begin{equation*}
\sum\limits_{j=1}^r3^{r-j}2^{j-1}=\sum\limits_{j=1}^r\binom{r}{j}(2^{j}-1).
\end{equation*}
It can be easily shown that both sums are equal to $3^r-2^r$ which completes the proof.\hfill{$\Box$}
\begin{theorem}
Let $X_0,X_1,X_2,X_3$ be independent random variables from the distribution $F\in\mathcal{F}$. If

\begin{equation}
 X_0 + X_{(2;3)}\overset{d}{=} X_{(3;3)}
\end{equation}
then $X\sim \mathcal{E}(\lambda), \lambda>0$.
\end{theorem}
\textbf{Proof.} Equalizing the respective densities we get:
\begin{eqnarray}
\nonumber \int\limits_{0}^{x}f(y)6F(x-y)(1-F(x-y))f(x-y)dy&=&3F^{2}(x)f(x)\\
\nonumber 6\int\limits_{0}^{x}f(y)F(x-y)(1-F(x-y))f(x-y)dy&=&6f(x)\int\limits_{0}^xF(y)f(y)dy\\
 \int\limits_{0}^{x}f(y)G(x-y)dy-\int\limits_{0}^{x}f(y)H(x-y)dy&=&f(x)\int\limits_{0}^xG(y)dy.\label{t2}
\end{eqnarray}

As in the previous theorem we shall prove by induction that from \eqref{t2} follows \eqref{fizvod}.

Differentiating both sides three times with respect to $x$ we get
\begin{eqnarray*}
 &&f''(x)G(0)+ f'(x)G'(0)+f(x)G''(0)+\int\limits_{0}^{x}f(x)G^{(3)}(x-y)dy\\
 &-&\Big(f''(x)H(0)+ f'(x)H'(0)+f(x)H''(0)+\int\limits_{0}^{x}f(x)H^{(3)}(x-y)dy\Big)\\
 &=& 3f''(x)G(0)+3f'(x)G'(0)+f(x)G''(0)+f^{(3)}\int\limits_{0}^xG(y)dy.
\end{eqnarray*}

Inserting $x=0$ we get

\begin{equation*}
 f'(0)G'(0)-f(0)H''(0)=3f'(0)G'(0),
\end{equation*}

and hence

\begin{equation*}
 f'(0)=-f^{2}(0).
\end{equation*}

Thus we got that \eqref{fizvod} is satisfied for $k=1$. Suppose now that \eqref{fizvod} is satisfied for $1\leq k\leq r-3$. Differentiating both sides of \eqref{t2} $r$ times we get

\begin{eqnarray*}
 &&\sum\limits_{j=0}^{r-1}f^{(r-1-j)}(x)G^{(j)}(0)+\int\limits_{0}^{x}f(y)G^{(r)}(x-y)dy\\&-& \sum\limits_{j=0}^{r-1}f^{(r-1-j)}(x)H^{(j)}(0)-\int\limits_{0}^{x}f(y)H^{(r)}(x-y)dy
 \\&=&\sum\limits_{j=1}^{r}\binom{r}{j}f^{(r-j)}(x)G^{(j-1)}(x)+f^{(r)}(x)\int\limits_{0}^xG(y)dy.
\end{eqnarray*}

Putting $x=0$ and eliminating zero terms we get

\begin{equation*}
 \sum\limits_{j=2}^{r}f^{(r-j)}(0)G^{(j-1)}(0)-\sum\limits_{j=2}^{r-1}f^{(r-1-j)}(0)H^{(j)}(0)=\sum\limits_{j=2}^{r}\binom{r}{j}f^{(r-j)}(0)G^{(j-1)}(0)
\end{equation*}

The terms for $j=r$ are equal on both sides so they cancel out. Extracting the summands containing $f^{(r-2)}(0)$ and grouping them on one side we get

\begin{eqnarray*}
&& f^{(r-2)}(0)G'(0)\bigg(1-\binom{r}{2}\bigg)\\&=&\sum\limits_{j=3}^{r-1}\bigg(\binom{r}{j}-1\bigg)f^{(r-j)}(0)G^{(j-1)}(0)+\sum\limits_{j=2}^{r-1}f^{(r-1-j)}(0)H^{(j)}(0)
\end{eqnarray*}

Using the induction hypothesis we obtain
%
\begin{eqnarray*}
 &&f^{(r-2)}(0)\bigg(1-\binom{r}{2}\bigg)\\&=&(-1)^{r-2}f^{r-1}(0)\left(\sum\limits_{j=3}^{r-1}(2^{j-1}-1)\bigg(\binom{r}{j}-1\bigg)
 -\sum\limits_{j=2}^{r-1}(3^{j}-2^{j+1}+1)\right),
\end{eqnarray*}

To prove that the condition \eqref{fizvod} is satisfied for $k=r-2$ it remains to show that
\begin{equation*}
 1-\binom{r}{2}=\bigg(-2+\sum\limits_{j=3}^{r-1}(2^{j-1}-1)\bigg(\binom{r}{j}-1\bigg)
 -(3^{j}-2^{j+1}+1)\bigg),
\end{equation*}

or, equivalently,

\begin{equation}
\sum\limits_{j=2}^{r-1}(2^{j-1}-1)\bigg(\binom{r}{j}-1\bigg)=\sum\limits_{j=2}^{r-1}(3^{j}-2^{j+1}+1)\label{identitet2}.
\end{equation}

It can easily be calculated that both sums are equal to $\frac{3^r}{2}-2^{r+1}+r+\frac{3}{2}$, which completes the proof.$\hfill \Box$

\begin{theorem}\label{druga}
Let $X_1,X_2,X_3,X_4$ be independent random variables from the distribution $F\in\mathcal{F}$. If

\begin{equation}
 X_{(2;3)}+\frac{1}{4}X_4\overset{d}{=} X_{(3;4)}
\end{equation}

then $X\sim \mathcal{E}(\lambda), \lambda>0$.

\end{theorem}
\textbf{Proof.}
Equalizing the respective densities we get:
\begin{eqnarray}
\nonumber \!\!\int\limits_{0}^{x}\!\!6F(x-y)(1\!-\!F(x-y))f(x-y)4f(4y)dy\!\!\!\!\!&=&\!\!\!\!\!12F^{2}(x)(1\!-\!F(x))f(x)\\
\nonumber 2\!\!\int\limits_{0}^{x}\!\!f(4y)F(x-y)(1\!-\!F(x-y))f(x-y)dy\!\!\!\!\!&=&\!\!\!\!\!f(x)\!\!\int\limits_{0}^x\!\!(F^2(y)\!-\!F^3(y))f(y)dy\\
 2\!\!\int\limits_{0}^{x}\!\!f(4y)(G(x-y)\!-\!H(x-y))dy\!\!\!\!\!&=&\!\!\!\!\!f(x)\!\!\int\limits_{0}^x\!\!(2G(y)\!-\!3H(y))dy.\label{t3}
\end{eqnarray}
As in the previous cases, we shall prove that from \eqref{t3} follows \eqref{fizvod}.

After evaluating the third derivative of both sides at zero and eliminating zero terms we get
\begin{equation*}
 8f'(0)G'(0)-2f(0)H''(0)= 6f'(0)G'(0)-3f(x)H''(0)
\end{equation*}

and hence again

\begin{equation*}
 f'(0)=-f^{2}(0).
\end{equation*}

The condition \eqref{fizvod} is thus valid for $k=1$. Suppose now that it is satisfied for $1\leq k\leq r-3$. Differentiating both sides of \eqref{t3} $r$ times we get
 \begin{eqnarray*}
&&2\sum\limits_{j=0}^{r-1}4^{r-1-j}f^{(r-1-j)}(4x)(G^{(j)}(0)-H^{(j)}(0))\\&+&2\int\limits_{0}^x4^{r}f^{(r)}(4y)(G(x-y)-H(x-y))dy\\
&=&\sum\limits_{j=1}^{r}\binom{r}{j}f^{(r-j)}(x)(2G^{(j-1)}(x)-3H^{(j-1)}(x))
+f^{r}(x)\int\limits_{0}^x(2G(y)-3H(y))dy.
\end{eqnarray*}
Letting $x=0$ and eliminating zero terms we get
\begin{eqnarray*}
&&2\sum\limits_{j=2}^{r-1}4^{r-j}f^{(r-j)}(0)G^{(j-1)}(0)-2\sum\limits_{j=2}^{r-1}4^{r-1-j}f^{(r-1-j)}(0)H^{(j)}(0)\\
&=&2\sum\limits_{j=2}^{r}\binom{r}{j}f^{(r-j)}(0)G^{(j-1)}(0)-3\sum\limits_{j=2}^{r-1}\binom{r}{j+1}f^{(r-1-j)}(0)H^{(j)}(0).
\end{eqnarray*}
The terms for $j=r$ in the first and the third sum coincide so they cancel out. Extracting terms containing $f^{(r-2)}(0)$ on one side we get 
\begin{eqnarray*}
&&2f^{(r-2)}(0)f^2(0)\Big(4^{r-2}-\binom{r}{2}\Big)\\&=&
2\sum\limits_{j=3}^{r-1}f^{(r-j)}(0)G^{(j-1)}(0)\Big(\binom{r}{j}-4^{r-j}\Big)\\&-&
\sum\limits_{j=2}^{r-1}f^{(r-1-j)}(0)H^{(j)}(0)\Big(3\binom{r}{j+1}-2\cdot4^{r-1-j}\Big).
\end{eqnarray*}
Applying the induction hypothesis to the right hand side we get 
\begin{eqnarray*}
&&2f^{(r-2)}(0)\Big(4^{r-2}-\binom{r}{2}\Big)\\&=&
(-1)^{r-2}f^{r-1}(0)\bigg(2\sum\limits_{j=3}^{r-1}(2^{j-1}-1)\Big(\binom{r}{j}-4^{r-j}\Big)\\&+&
\sum\limits_{j=2}^{r-1}(3^j-2^{j+1}+1)\Big(3\binom{r}{j+1}-2 \cdot4^{r-1-j}\Big)\bigg).
\end{eqnarray*}
It remains to prove
\begin{equation}
\sum\limits_{j=2}^{r-1}(2^{j}-2)\bigg(4^{r-j}-\binom{r}{j}\bigg)=\sum\limits_{j=2}^{r-1}(3^{j}-2^{j+1}+1)\bigg(3\binom{r}{j+1}-2\cdot4^{r-j-1}\bigg)\label{identitet3}.
\end{equation}

It can easily be calculated that both sums are equal to $\frac{4^r}{3}-3^{r}+2^{r}-\frac{1}{3}$, which completes the proof.$\hfill \Box$
\section{Discussion}

In this paper we stated and proved three characterizations theorems for exponential distribution. Although they have similar formulation they could be considered as different types.
We suppose that these three types of characterizations could be valid for any order statistics. For example, the theorem \ref{prva} could be generalized that $\sum_{j=1}^k\frac{X_{j}}{n-j+1}$ and $X_{(k;n)}$ are 
identically distributed. Other two characterizations could also be generalized in a similar manner. So far, the generalization of the theorem \ref{druga} for the special case of $k=n$ (consecutive maxima)
has been considered  in \cite{chakraborty}. We hope to extend the results in the future.

\end{document}